\documentclass[10pt]{article}
\usepackage{amsmath,amsfonts,amsthm,amssymb,amscd}
\numberwithin{equation}{section}

\newcommand{\be}{\begin{equation}}
\newcommand{\ee}{\end{equation}}
\newcommand{\bs}{\begin{split}}
\newcommand{\es}{\end{split}}
\newcommand{\ba}{\begin{align}}
\newcommand{\ea}{\end{align}}
\newcommand{\basl}[1]{\begin{align}\begin{split}\label{#1}}
\newcommand{\bas}{\begin{align}\begin{split}}

\newtheorem{theo}{Theorem}[section]
\newtheorem{prop}[theo]{Proposition}
\newtheorem{lemm}[theo]{Lemma}

\newcommand\fpr{\hfill$\Box$\null}

\newcommand\R{\mathbb{R}}
\newcommand\C{\mathbb{C}}

\title{Integral formulas for the Weyl and anti-Wick symbols }

\author{ L. Amour and  J. Nourrigat}

\date{Universit\'e de Reims, France}

\begin{document}

\maketitle

\begin{abstract}
\noindent
The first purpose of this article is to provide conditions for a bounded operator in $L^2(\R^n)$ to be the Weyl (resp. anti-Wick) quantization of a bounded continuous symbol on $\R^{2n}$. Then, explicit formulas for the Weyl (resp. anti-Wick) symbol are proved.  Secondly, other formulas for the Weyl and anti-Wick symbols involving a kind of Campbell Hausdorff formula are obtained. A point here is that these conditions and explicit formulas depend on the dimension $n$ only through a Gaussian measure on $\R^{2n}$ of variance $1/2$ in the Weyl case (resp. variance $1$ in the anti-Wick case) suggesting that the infinite dimension setting for these issues could be considered.  Besides, these conditions are related to iterated commutators recovering in particular the Beals characterization Theorem.
\end{abstract}

\parindent=0pt

\


\


\parindent=0pt
\parindent = 0 cm

\parskip 10pt
\baselineskip 12pt

 \section{Statement of the results.}\label{s1}

Our purpose in this paper is to answer the following issues when $A$ is a  bounded  operator in 
 $L^2 (\R^n)$:

 1. Under what condition,  the Weyl symbol of  $A$ is a bounded continuous function on $\R^{2n}$ ?

 2. Under what condition,   $A$ is the anti-Wick quantization of a bounded continuous function on $\R^{2n}$ ?

 3. Can one express the Weyl and anti-Wick symbols with a  Campbell-Hausdorff type formula ?
 
These questions are also mentioned in  \cite{K-R}. A positive answer to the first two questions (a sufficient condition) is given by Theorem \ref{t1}  for the Weyl symbol and by Theorem \ref{t3}
for the anti-Wick symbol. Concerning Question 3, an answer is given combining  Theorem \ref{t1} or Theorem \ref{t3}  with Theorem  \ref{transfo-beals}. For instance, one gets (\ref{AWCH})
in the case of  the anti-Wick quantization.

Let us emphasize that it is also our aim to provide conditions and formulas for the two symbols where the dimension is not explicitly written down using Gaussian measures. Thus, the dimension can tends to infinity and more importantly, one can expect to consider the infinite dimension setting for theses purposes in a further work. 

In the article \cite{HS-RD}, H\"ubner and Spohn have introduced the so-called scattering identification operator (unbounded operator) and considered the wave operator in the context of the infinite dimension, concerning some physical model. With L. Jager, we have proved in
 \cite{A-J-N} that the first operator (scattering identification) can also be defined with the anti-Wick quantization. That is, this operator is defined by an integral formula and it would be interesting to know if the second operator (wave) shares that property. One may expect that the present work, adjusted to the infinite dimension case, could  give an answer to that question.

We begin by recalling the two quantizations under consideration here. See, $e.g.$, 
 Berezin \cite{FAB},  Combescure and Robert \cite{C-R}, Folland \cite{F}, H\"ormander \cite{HO}, 
 Lerner \cite{Ler}, Unterberger \cite{U}\cite{Unt-OH}, $\dots$ for these purposes.

One can define the Weyl symbol of a continuous operator $A$ from  ${\cal S}(\R^n)$
to ${\cal S}'(\R^n)$ as an element of  ${\cal S}'(\R^{2n})$. More precisely, denoting by $K_A(x , y)$ the distributional kernel 
of $A$ (belonging to ${\cal S}'(\R^{2n}))$,
the  Weyl symbol of $A$ is defined by:
%------
\be\label{noyau-symbole}  F(x , \xi) =   \int _{\R^n} e^{ - i v\cdot \xi} K_A ( x + v/2,  x - v/2)
dv \ee
%----
where the above integral is actually a Fourier transform in sense of tempered distributions.

The second kind of quantization under consideration in this paper is the anti-Wick one, which involves  coherent states.
These states are a standard family of functions  $\Psi_X$ indexed by $X = (x , \xi)\in \R^{2n}$ defined by:
 %-----
 $$ \Psi _{x , \xi} (u) =  \pi ^{-n/4} e^{ - \frac {1} {2} |u-x|^2
 + i u \cdot \xi -  \frac {i} {2} x \cdot \xi}. $$
 %------
One has:
 %----
 \be\label{PSEC}  < \Psi_X , \Psi_Z >
 = e^{ -  \frac {|X- Z|^2} {4} + \frac {i} {2}
\sigma (X , Z) } \ee
%----------
where  $\sigma $ is the symplectic form, $\sigma ( (x , \xi), (y, \eta)) = y \cdot \xi - x \cdot \eta$. 
For any function $G$,  bounded and continuous on  $\R^{2n}$, one defines
 the anti-Wick operator $Op^{AW}(G)$ associated with the function $G$, as the unique operator satisfying for all $f$ and $g$ in $L^2(\R^n)$:
 %-----
 \be\label{def-AW}  <Op^{AW}(G) f , g> = (2\pi )^{-n} \int _{\R^{2n}} G(Z)
 < f, \Psi_Z> < \Psi _Z , g> dZ. \ee
 %-------

Before developing Question 3 with Campbell Hausdorff type formulas and iterated commutators, we present answers to the first two questions addressed above, in the spirit of the characterization of A. Unterberger \cite{U}.

\begin{theo}\label{t1}
Let $A$ be a bounded operator in $L^2(\R^n)$ and assume that:
 %---
 \be\label{hyp-Weyl}   \sup _{X\in \R^{2n}}  \pi ^{-n} \int _{\R^{2n} }
 \left |\frac {   <  A \Psi_{X+Z}  , \Psi_{X-Z}  >} {  <   \Psi_{X+Z}  , \Psi_{X-Z}  >} \right |
 e^{ - |Z|^2} dZ <  \infty. \ee
 %------
 Then, the Weyl symbol $F$ of $A$ is a bounded continuous function on $\R^{2n}$  given by the formula:
 %----
\be\label{symbole-weyl}  F(X) =  \pi ^{-n}  \int _{\R^{2n}}
\frac {   <  A \Psi_{X+Z}  , \Psi_{X-Z}  >} {  <   \Psi_{X+Z}  , \Psi_{X-Z}  >}
  e^{ -|Z|^2}  dZ. \ee
%-----

\end{theo}

We underline that hypothesis (\ref{hyp-Weyl}) and equality (\ref{symbole-weyl}) are uniform with respect to the dimension since these expressions only use the variance $1/2$  Gaussian measure which also exists in infinite dimension. For a reader interested only in the finite dimensional case, one can check that, using (\ref{PSEC}),
the hypothesis  (\ref{hyp-Weyl}) can also be written as:
%---
$$  \sup _{X\in \R^{2n}} \int _{\R^{2n} }
 |< A \Psi_{X+Z} ,  \Psi_{X-Z}>| dZ <  \infty $$
%----
and (\ref{symbole-weyl}) can also be expressed as:
%---
$$ F(X)  =  \pi ^{-n}  \int _{\R^{2n}}  <  A \Psi_{X+Z}  , \Psi_{X-Z}  > e^{i \sigma (X , Z) }
dZ.$$
%----

Next we turn to conditions ensuring the existence of a bounded continuous symbol with the anti-Wick quantization for a given bounded operator. 

Firstly, we give a necessary condition.

 \begin{theo}\label{t2}  If $A= Op^{AW}(G)$ with  $G$ bounded and continuous on $\R^{2n}$ then:
 %----
 \be\label{CN}  |<  A \Psi_X  , \Psi_{Y}  > | \leq \Vert G \Vert _{\infty } e^{ - \frac {1} {8}
| X -  Y |^2 }.\ee
%--------

 \end{theo}

In particular, the following estimate holds,
 %---
 $$   | < A \Psi_{X+Y}   , \Psi_{X -Y}  >| \leq \Vert G \Vert _{\infty }
 e^{ - \frac {1} {2} |  Y |^2 }$$
%--------
for all $X$ and $Y$ in $\R^{2n}$, when $A= Op^{AW}(G)$ where  $G$ is bounded and continuous on $\R^{2n}$.

Secondly, we provide a sufficient condition, together with an integral expression for a possible anti-Wick symbol.

 \begin{theo}\label{t3} Suppose that  $A$ is a bounded operator in $L^2 (\R^n)$ satisfying:
 %----
 \be\label{CSAW}  \sup  _{X\in \R^{2n}} (2\pi )^{-n} \int _{\R^{2n}}
  \left |\frac {   <  A \Psi_{X+Z}  , \Psi_{X-Z}  >} {  <   \Psi_{X+Z}  , \Psi_{X-Z}  >} \right |
 e^{ - \frac {|Z|^2} {2} } dZ < \infty. \ee
 %-------
 Then, there exists a bounded continuous function $G$ on $\R^{2n}$ such that
 $A = Op^{AW} (G)$.  This function is given by:
   %----
   \be\label{formule-symb-AW}
   G(X) =  (2\pi )^{-n}  \int _{\R^{2n}}
   \frac {   <  A \Psi_{X+Z}  , \Psi_{X-Z}  >} {  <   \Psi_{X+Z}  , \Psi_{X-Z}  >}
    \ e^{- \frac {1} {2} |  Z |^2 } dZ.\ee
   %----
\end{theo}

Note that condition (\ref{CSAW}) can also be written as:
 %----
 $$ \sup  _{X\in \R^{2n}} (2\pi )^{-n} \int _{\R^{2n}}
   | < A \Psi_{X+Z}   , \Psi_{X -Z}  >|\   e^{  \frac {1} {2} |  Z |^2 }
   dZ < \infty$$
%------
and equality (\ref{formule-symb-AW}) is also:
%---
$$  G(X) =  (2\pi )^{-n}  \int _{\R^{2n}}
    < A \Psi_{X+Z}   , \Psi_{X -Z}  > e^{ i\sigma (X , Z)}
   \ e^{  \frac {1} {2} |  Z |^2 } dZ. $$
%---
 
We now give answers to Question 3 in the spirit of Beals characterization \cite{Bea}.

We first note that the function appearing in  (\ref{symbole-weyl})  and (\ref{formule-symb-AW})
can also be written with an expression comparable to the  Campbell-Hausdorff formula.  To this end, we denote by $ \Phi_S(Z)$ the differential operator:
%-----
\be\label{PhiS}  \Phi_S(Z) = \sum _{j = 1}^n \left ( z_j u_j + \zeta_j \frac {1} {i}
\frac {\partial } {\partial u_j}\right ) \ee 
%---
for all $Z = (z , \zeta)\in \R^{2n}$.
\begin{theo}\label{transfo-beals} For all $X$ and $Z$ in $\R^{2n}$, we have the following equality:
%----
\be\label{it-comm}  \frac {   <  A \Psi_{X+Z}  , \Psi_{X-Z}  >} {  <   \Psi_{X+Z}  , \Psi_{X-Z}  >}
= \Big < e^{- \Phi_S(Z) }  A e^{\Phi_S(Z) }   \Psi_{X},
 \Psi_{X}   \Big > . \ee
%------

\end{theo}

One notices that the operator  $e^{\Phi_S(Z) }$  is unbounded in $L^2(\R^n)$. Nevertheless, its action on coherent states is well defined. It is  given by:
%---
\be\label{segal-EC}  e^{\Phi_S(Z) } \Psi_X = e^{ \frac{1} {2} |Z|^2 + Z \cdot X -  \frac{i} {2} \sigma (Z, X)}
     \Psi_{X+Z}.\ee
%------
One observes that (\ref{it-comm}) only uses iterated commutations according to the Campbell Hausdorff formula.
One has in  the sense of formal series: 
%-----
$$ e^{-\Phi_S(Z) }  A e^{\Phi_S(Z) }   = \sum _{m= 0}^{\infty }  \frac{(-1)^m } {m!}
   \big ( {\rm ad}\, \Phi_S(Z) \big ) ^m A. $$
%---
Inserting this equality in  (\ref{symbole-weyl})  and in (\ref{formule-symb-AW}), one then obtains an expression for the Weyl and anti-Wick symbol  of $A$ with an integral using only iterated commutations. 
To be specific, for the possible anti-Wick symbol $G$ of an operator $A$, if the integral is absolutely converging, one has:%-----
\be\label{AWCH}   G(X) =  (2\pi )^{-n}  \int _{\R^{2n}} 
\Big < e^{- \Phi_S(Z) }  A e^{\Phi_S(Z) }   \Psi_{X},
 \Psi_{X}   \Big >  \ e^{  - \frac {1} {2} |  Z |^2 } dZ. \ee 
%----

One can also give a result (see below) specific to the finite dimension, close to the Beals characterization theorem \cite{Bea}. 
 Denote by $V_1,\dots,V_n$ the canonical basis of
$\R^n$ and set  $a(V_j) = x_j + \partial _{x_j}$. For all multi-indexes
$\alpha$ and any bounded operator $A$ from ${\cal S}(\R^n)$
to ${\cal S}'(\R^n)$, set $({\rm ad} \,a(V)) ^{\alpha} A =
{\rm ad}\,  a(V_1)^{\alpha _1} \dots {\rm ad}\,  a(V_n)^{\alpha _n}A$.

\begin{prop}\label{beals}  For any operator $A$ from ${\cal S}(\R^n)$
to ${\cal S}'(\R^n)$ with iterated commutators 
$({\rm ad}\, a(V)) ^{\alpha} A$ bounded 
in $L^2 (\R^n)$ when $|\alpha | \leq 2n+1$, one has:
%------
\be\label{norme-beals}    \int _{\R^{2n}}
  |<  A \Psi_{X+Z}  , \Psi_{X-Z}  >|  dZ \leq C_n
  \sum _{ |\alpha | \leq 2n+1 } \Vert ({\rm ad}\, a(V)) ^{\alpha} A \Vert \ee
%---
where the constant  $C_n$ depends on the dimension  $n$ but not on the operator $A$.

\end{prop} 
We then recover in particular a well known result of  R. Beals \cite{Bea} 
for the Weyl symbol $F$ of an operator $A$:
%---
$$ \sup _{X\in \R^{2n}} |F(X)| \leq C_n  \sum _{ |\alpha | \leq 2n+1 } 
\Vert ({\rm ad}\,  a(V)) ^{\alpha} A \Vert. $$
%----- 
We also mention that there is a result of  C. Rondeaux \cite{R} where above, the supremum bound is replaced by the $L^1$ norm and the operator norm is replaced by the trace class norm.

Theorem  \ref{t1} and  Theorem \ref{t3} rely on a inversion result for the heat operator (Theorem \ref{t-2-1}) that is proved in Section \ref{s2}. On the basis of this result (Theorem \ref{t-2-1}),  Theorems
\ref{t1} and  \ref{t3}  are proved in Sections  \ref{s3} and  \ref{s4}. 
 Equality (\ref{segal-EC}),   Theorem \ref{transfo-beals} 
and Proposition \ref{beals} are proved in Section \ref{s5}.

 \section{Inversion of the heat operator.}\label{s2}

For each $\lambda >0$,  the heat operator $H _{\lambda }$ is defined for all bounded continuous function $F$ on $\R^{2n}$ by:
%-----
$$ (H _{\lambda } F) (z , \zeta) =   (2 \pi \lambda ) ^{-n}  \int _{\R^{2n}} F(x, \xi )
e^{ -  \frac {1} {2 \lambda } \big (  (x-  z)^2 +  (\xi  - \zeta )^2 \big )  } dx d\xi
\hskip 2cm (z , \zeta ) \in \C^{2n} $$
%------
where  $z^2  = z_1 ^2 + \cdots + z_n ^2$ for all $z = ( z_1 , \cdots z_n)\in\C^n$. 
This function is holomorphic on $\C^{2n}$.

The purpose of Theorem \ref{t-2-1} is to prove that the image of the operator $H _{\lambda }$ 
contains some specific space playing a role in the following sections. We also use another function $ S _{\lambda } F$
on $\C^{2n}$ defined by:
%----
\be\label{rotat}  (S _{\lambda } F) (z , \zeta) =  (H _{\lambda } F) (i \zeta, -iz ).\ee 
%-----
That is:
%---
\be\label{S-lambda-1} (S _{\lambda } F) (z , \zeta) =   (2 \pi \lambda ) ^{-n}  \int _{\R^{2n}} F(x, \xi )
e^{ -  \frac {1} {2 \lambda } \big (  (x - i \zeta)^2 +  (\xi  + i z)^2 \big )  } dx d\xi
\hskip 2cm (z , \zeta ) \in \R^{2n}.  \ee
 %------

\begin{theo}\label{t-2-1}   Let $\Phi$  be a holomorphic function on  $\C^{2n}$. Set $\lambda >0$ and suppose that:
%----
\be\label{S-lambda-2} \sup _{(x , \xi) \in \R^{2n}  } (2\pi \lambda ) ^{-n} \int _{\R^{2n}}
|\Phi ( x + i \xi + z + i \zeta ,   x - i \xi - z + i \zeta ) | e^{- \frac {|z|^2 + |\zeta |^2  } {2\lambda }} dz d\zeta  < \infty. \ee
%-------
Define the function $F$ on  $\R^{2n}$  by:
%-----
\be\label{S-lambda-3} F(x , \xi) =  (2\pi \lambda ) ^{-n} \int _{\R^{2n}}
\Phi ( x + i \xi + z + i \zeta ,   x - i \xi - z + i \zeta )  e^{- \frac {|z|^2 + |\zeta |^2  } {2\lambda }} dz d\zeta.
\ee
 %-----
 Then,  $F$ is a bounded continuous function on $\R^{2n}$ satisfying:
  %----
 \be\label{H-lambda} (H_{\lambda } F ) (z , \zeta ) = \Phi ( z + i \zeta  ,
  z - i \zeta ) \hskip 2cm (z , \zeta ) \in \C^{2n} \ee
 %------
 %----
 \be\label{S-lambda-4} (S_{\lambda } F ) (z , \zeta ) = \Phi ( z + i \zeta  ,
 - z + i \zeta ) \hskip 2cm (z , \zeta ) \in \C^{2n}. \ee
 %------
 \end{theo}

 {\it Proof.} The fact that $F$ is well defined, bounded and continuous, is a direct consequence of hypothesis (\ref{S-lambda-2}).
 For all $(x , \xi)\in\R^{2n}$,  let $\Psi_{(x , \xi)} $ be the holomorphic function on $\C^{2n}$ defined by:
 %----
 $$ \Psi _{(x , \xi)}  (Z_1 , Z_2) =  \Phi  (Z_1 , Z_2) e^{   \frac {1} {2 \lambda }
 \big ( Z_1 Z_2 + x^2 + \xi^2 - x ( Z_1 + Z_2) + i \xi ( Z_1 - Z_2) \big ) }.$$
 %---
 Set $\varphi$  the mapping from $\R^{2n}$ to $\C^{2n}$ defined by
 $\varphi (s, t) = ( s+i t, -s + i t)$.
Equality  (\ref{S-lambda-3}) can then be written as:
 %----
 $$ F(x , \xi) =  (2\pi \lambda ) ^{-n}  \int_{\R^{2n}}
   \Psi _{(x , \xi)}  \big (  ( x+ i \xi , x - i \xi) + \varphi (s, t) \big )
 ds dt. $$
 %----
 We then apply Lemma \ref{DSS} below with $n$ replaced by
 $2n$, $E = \{ (x + i \xi, x - i \xi), \ (x , \xi) \in \R^{2n} \}$, together with the above functions  $\varphi$ and $\Psi =  \Psi _{(x , \xi)}$. According to hypothesis (\ref{S-lambda-2}), for all compact sets $K$ 
 of $\R^{2n}$, with $(x, \xi)$ being fixed:
 %----
 $$ \sup _{(y, \eta)\in K}  \int _{\R^{2n}}
|\Psi _{(x , \xi)} \big ( (y+ i \eta ,   y - i \eta) + \varphi (s , t) \big ) | ds dt < \infty.$$
%----
 Consequently, hypothesis  (\ref{hyp-lemm})  of Lemma \ref{DSS} is satisfied. 
 According to that Lemma, one has:
 %----
 $$ F(x , \xi) =  (2\pi \lambda ) ^{-n}  \int _{\R^{2n}}
 \Psi _{(x , \xi)}   \big ( \varphi (s, t) \big )
 ds dt.$$
 %----
 This equality can be written as:
 %----
 $$ F(x , \xi) =   (2\pi \lambda ) ^{-n} \int _{\R^{2n}}
 \Phi ( z + i \zeta ,   - z + i \zeta )
 e^{ -  \frac {1} {2 \lambda } \big (  (z - i \xi)^2 +  (\zeta   + i x)^2 \big )  }
 dz d\zeta.  $$
 %-------
Set:
 %----
 $$ G(z, \zeta) =   \Phi ( z + i \zeta ,   - z + i \zeta )
 e^{- \frac {|z|^2 + |\zeta |^2  } {2\lambda }}
 \hskip 2cm
 H(x , \xi) = F( \xi, -x)\  e^{- \frac {|x|^2 + |\xi |^2  } {2\lambda }}.  $$
%----
Hypothesis (\ref{S-lambda-2}) with  $(x , \xi)= (0, 0)$ shows that the function $G$ belongs to $L^1(\R^{2n})$.
The above equality reads as: 
%----
$$  H(x , \xi) = ( 2\pi \lambda ) ^{-n} \int _{\R^{2n}} G(z, \zeta)
e^{- \frac {i} {\lambda }  ( x\cdot z + \xi \cdot \zeta)} dz d \zeta.$$
%--------
Thus, $H$ is the Fourier transform depending on the parameter
$\lambda$ of $G$. Since $F$ is bounded and continuous then $H\in L^1(\R^ {2n})$.
Hence,
%----
$$  G(x , \xi) = ( 2\pi \lambda ) ^{-n} \int _{\R^{2n}} H(z, \zeta)
e^{ \frac {i} {\lambda }  ( x\cdot z + \xi \cdot \zeta)} dz d \zeta.$$
%--------
This equality is equivalent to  (\ref{S-lambda-4}) for real  $(z , \zeta)$.
Equality (\ref{S-lambda-4}) then follows for $(z , \zeta)\in \C^{2n}$ using the holomorphic properties of the two hand sides. Equality (\ref{H-lambda}) then holds according to (\ref{rotat}).  \fpr

The Lemma below would be in the particular case of the dimension one, the result for changing integration contours with holomorphic functions when these integration contours are parallel lines.

\begin{lemm}\label{DSS}  Set $E$ a $n-$dimensional real subspace of $\C^n$.
Let $\varphi $ be a linear map from $\R^n$ to $\C^n$ satisfying 
 ${\rm det } \varphi ' \not=0$. Assume that $E \oplus {\rm Im} \varphi = \C^n$. 
Set $\Psi$ a holomorphic function on $\C^n$.
Suppose that, for any compact set $K$ of $E$:
%----
\be\label{hyp-lemm}  \sup _{a\in K} \int _{\R^n} | \Psi (a + \varphi (t ))|  d t < \infty.  \ee 
%---
Set:
%---
$$ I(a) = \int _{\R^n} \Psi (a + \varphi (t )) d t $$
%--
for any $a\in E$ where $dt $ is the  Lebesgue measure on $\R^n$. Then $I$ is independent of  $a\in E$.

\end{lemm}

{\it Proof of the Lemma.} For all $R>0$, set:
%---
$$ K_R = \{ t \in \R^n , \ \ |t_j| \leq R , \ \ \ \ 1 \leq j \leq n \}.$$
%------
Set $\partial K_R$ the boundary of $K_R$ with the canonical measure $d \mu_R$. 
For all $k\leq n$, we have:
%----
$$ \left | \int _{K_R }  \frac { \partial } {\partial t_k} \Psi (a + \varphi (t )) dt 
\right  |  \leq \int _{\partial K_R } |\Psi (a + \varphi (t )) | d \mu_R(t).  $$
%-----
Since $\Psi$ is holomorphic:
%---
$$ \left | \sum _{j= 1}^n  \frac { \partial \varphi_j  } {\partial t_k}
  \int _{K_R}  (\partial _j \Psi) (a + \varphi (t )) dt \right  | 
  \leq \int _{\partial K_R } |\Psi (a + \varphi (t )) | d \mu_R(t).   $$
  %--
As the determinant of the matrix $ \frac { \partial \varphi_j  } {\partial t_k} $
is non vanishing, one gets for all  $j \leq n$:
  %----
  $$ \left |   \int _{K_R}  (\partial _j \Psi) (a + \varphi (t )) dt  \right  |  
   \leq C \int _{\partial K_R } |\Psi (a + \varphi (t )) | d \mu_R(t)  $$
  %--
with $C>0$. Using that  $\Psi$ is holomorphic, there exists another constant $C>0$ such that,  for all $Z\in \C^n$: 
   %---
  \be\label{moyenne} |\Psi (Z)| \leq C \int _{B(Z, 1)} |\Psi (W)| dW. \ee 
  %-------
 Consequently, with another constant  $C$ and another $\rho >0$:
  %-------
  $$ \int _{\partial K_R } |\Psi (a + \varphi (t )) | d \mu_R(t) 
  \leq C \int _{B_E(a, \rho) \times F ( \partial \Gamma _R , \rho) } 
  |\Psi (a'  + \varphi (t )) | \  d\mu_E (a') dt $$
  %----
  where $B_E(a, \rho) $ is the ball of $E$  centered at $a$ with radius $\rho$, $\mu_E$ is the canonical measure on $E$ and  $F ( \partial \Gamma _R , \rho)$ 
  is the set of all points in  $\R^n$ whose distance to $\partial \Gamma _R $ 
  is smaller than $\rho$. One then deduces that, for all  $a$ and $b$ in $E$, for any $R>0$:
  %----
  $$  \left | \int _{K_R}  \big (  \Psi (a + \varphi (t )) - \Psi (b + \varphi (t )) 
  \big ) dt   \right  |  \leq  C  \int _{F_E([a, b], \rho) \times F ( \partial \Gamma _R , \rho) }
  |\Psi (a'  + \varphi (t )) | \  d\mu_E (a') dt   $$
  %---
  where $F_E([a, b], \rho)$ is the set of all points in  $E$ whose distance to the interval 
  $[a, b]$ is smaller than $\rho$. 
  The right hand side tends to zero (the constants $C$ and $\rho$ are independent of $R$) as $R$ goes to infinity.
 One then concludes that 
  $I(a) = I(b)$ which proves the Lemma. 
  \fpr

Let us also mention another property of the operator $S_{\lambda}$. Denote by  $M_{\lambda}$ the following multiplication operator:
%-----
$$ (M_{\lambda}F) (x , \xi) = e^{-  \frac {1} {2 \lambda } ( |x|^2 +  |\xi |^2) }
F(x , \xi). $$
%----
One then sees that $M_{\lambda}\ S_{\lambda} F = T_{\lambda}\ M_{\lambda} F$ where:
%----
$$  (T_{\lambda}F) (z , \zeta) = ( 2\pi \lambda ) ^{-n} \int _{\R^{2n}}
F(x , \xi) \  e^{  \frac {i} { \lambda } ( z \cdot \xi - x \cdot \zeta ) }
dx d\xi.$$
%-----
One notes that $T_{\lambda}$ is an isometry from $L^2(\R^{2n} )$ into itself,
whose square is the identity operator. It is the symplectic Fourier transform with parameter (see \cite{C-R}). Hence,  $S_{\lambda}$ is an isometry from
$E = \{ F , M_{\lambda}F\in L^2(\R^{2n} ) \} $ into itself, whose square is the identity. These properties are not further used in this work.

 \section{Weyl symbol.}\label{s3}

It is knwon (\cite{F}, page 139)  that, for any bounded operator $A$ from  ${\cal S}(\R^n)$ to ${\cal S}'(\R^n)$,
there exists an homomorphic function   $B(A)$ on $\C^{2n}$  such that, for all  $X = (x , \xi)$ and $Y = (y , \eta)$  identified with elements of $\C^n$:
%----
\be\label{bisymbole}  \frac {   <  A \Psi_X  , \Psi_{Y}  >} {  <   \Psi_X  , \Psi_{Y}  >} =
B(A) ( x+ i \xi, y - i \eta). \ee
%----
Set $A$ a bounded operator in
$L^2(\R^n)$ satisfying (\ref{hyp-Weyl}).  Assumption (\ref{hyp-Weyl}) then reads as:
 %----
 $$  \sup _{X\in \R^{2n}}  \pi ^{-n} \int _{\R^{2n} }
 \left | B(A) ( X+Z , \overline X - \overline Z ) \right |
 e^{ - |Z|^2} dZ <  \infty. $$
 %-------
That is, the function  $\Phi = B(A)$ satisfies the assumption of Theorem 
 \ref{t-2-1} with $\lambda = 1/2$. According to that Theorem, the function $F$ defined by  (\ref{S-lambda-3}), that is to say by (\ref{symbole-weyl}) when $\lambda = 1/2$ and $\Phi = B(A)$, is continuous and bounded on  $\R^{2n}$ and verifies (\ref{S-lambda-4}) with $\lambda = 1/2$.
Namely:
 %----
 \be\label{egal-bisymb-1}  (S_{1/2}F) (Z) =   B( A) (Z , - \overline Z).    \ee 
 %---
Inverting the  transform in
 (\ref{noyau-symbole}), one associates with   $F$ some tempered distribution on  
 $\R^{2n}$ which is the kernel $K_{A'} $ of some operator $A'$, bounded from ${\cal S}(\R^n)$ 
to ${\cal S}'(\R^n)$,  whose Weyl symbol is $F$.  
In view of the Lemma below, one has for all $Z\in\C^n$:
 %----
 \be\label{egal-bisymb-2}  B( A') (Z , - \overline Z) = (S_{1/2}F) (Z).   \ee
 %----
We recall Proposition 1.69 in \cite{F}.  If  $G$ is a holomorphic function on  $\C^n \times \C^n$
and if  $G (Z, - \overline Z)=0$ for all $Z\in\C^n$ then  $G$ is entirely vanishing. 
Indeed, all the  $\partial ^{\alpha }G$ are vanishing at the origin by iteration.  
Thus, from (\ref{egal-bisymb-1}) and (\ref{egal-bisymb-2}), 
      $ B( A - A')$ is identically zero. Then,
 $ < (A-A')\Psi_X , \Psi_Y >= 0$ for all $X$ and $Y$ in $\R^{2n}$ proving that
 $A-A'= 0$. Therefore $F$ is the Weyl symbol of  $A$.

 \begin{lemm} Set $F$ a bounded continuous function on  $\R^{2n}$. Let
 $A' $ be the bounded operator from    ${\cal S}(\R^n)$ to ${\cal S}'(\R^n)$
  with  $F$ as Weyl symbol.  Then, (\ref{egal-bisymb-2}) holds where
  $S_{1/2}$  is defined in (\ref{S-lambda-1}). 

  \end{lemm}

 {\it Proof of the Lemma.}  In view of (\ref{noyau-symbole}),  the operator $A'$, bounded from  ${\cal S}(\R^n)$ to ${\cal S}'(\R^n)$)
 with a Weyl symbol equal to  $F$, satisfies for all  $X$ and $Y$ in $\R^{2n}$:
 %-----
 \be\label{Weyl-A} < A' \Psi_X  , \Psi_Y  >  =
  \pi ^{-n}  \int _{\R^{2n}} F(Z) <\Sigma _Z \Psi_X  , \Psi_Y  >  dZ \ee 
 %----
 where $ \Sigma _Z$ for each $Z = (z, \zeta)\in\R^{2n}$ is the operator acting in 
 $L^2(\R^n)$ defined by:
 %----
 $$ \Sigma _Z f (u) = e^{2i ( x - z)\cdot \zeta} f ( 2z - x)
 \hskip 2cm    f\in L^2(\R^n). $$
 %------
 In particular:
 %---
 $$ \Sigma _Z \Psi_X = e^{i \sigma (X , Z)} \Psi_{ 2Z - X} $$
 %------
 and therefore, using (\ref{PSEC}):
 %----
 \be\label{Weyl-B}  < \Sigma _Z \Psi_X  ,\Psi_{-X} >  
 = e^{2 i \sigma (X , Z)-  |Z|^2 }.  \ee 
 %---
One deduces from  (\ref{bisymbole}), (\ref{PSEC}), (\ref{Weyl-A}), (\ref{Weyl-B}), 
and (\ref{S-lambda-1}):
 %----
 \begin{align*}
 B(A') (Z , - \overline X )  &=  \frac {   <  A' \Psi_X  , \Psi_{-X}  >} {  <   \Psi_X  , \Psi_{-X}  >}  \\ 
 &=  e^{ |X|^2} <A'  \Psi_X ,  \Psi_{-X} > \\ 
 %----
&=  \pi ^{-n} e^{ |X|^2}   \int _{\R^{2n}} F(Z) <\Sigma _Z \Psi_X  , \Psi_{-X}   >  dZ  \\ 
 %------
 &= \pi ^{-n}  \int _{\R^{2n}} F(Z)  e^{ |X|^2 + 2 i \sigma (X , Z)-  |Z|^2 }  dZ    = S_{1/2} F (X)
 %-----
 \end{align*}
which proves the Lemma. \fpr

  \section{Anti-Wick symbol.}\label{s4}

{\it Proof of Theorem \ref{t2}. } If $A = Op^{AW}(G)$ then, according to (\ref{def-AW}):
%----
$$ <  A \Psi_X  , \Psi_{Y}  > = (2\pi )^{-n}  \int _{\R^{2n}}  G(Z)
 < \Psi_X , \Psi_Z> < \Psi _Z , \Psi_{Y} > dZ. $$
%---
Using (\ref{PSEC}):
%----
\be\label{egalite}   < A \Psi_X  , \Psi_{Y}  >  = (2\pi )^{-n}  \int _{\R^{2n}}  G(Z)
e^{ - \frac {1} {4} \big ( |X - Z|^2 + |Y - Z|^2 \big ) + \frac {i} {2}
\sigma (X - Y , Z)       }dZ.\ee
%------
From the parallelogram identity:
%----
$$ | < A \Psi_X  , \Psi_{Y}  >| \leq (2\pi )^{-n}  \int _{\R^{2n}}  |G(Z)|
e^{ - \frac {1} {2} \left  | Z - \frac {X+Y} {2} \right  |^2 - \frac {1} {8}|X-Y|^2 } dZ.$$
%----
One then obtains (\ref{CN}).\fpr

{\it Proof of Theorem \ref{t3}. } Set   a bounded operator $A$ in
$L^2(\R^n)$ satisfying (\ref{CSAW}). That is, the function $B(A)$ defined by (\ref{bisymbole}) satisfies the assumption of Theorem  \ref{t-2-1} with $\lambda = 1$. According to that Theorem, the function $G$ defined by
 (\ref{S-lambda-3}) with $\lambda = 1$ and $\Phi = B(A)$, or equivalently by
  (\ref{formule-symb-AW}),  is bounded and continuous on $\R^{2n}$ and satisfies (\ref{S-lambda-4})
with $\lambda = 1$. In other words:
 %-----
 $$ (S_{1}G) (Z)  =  B(A) (Z , - \overline Z).$$
 %-------
One can define an operator $A' = Op^{AW} (G)$ bounded in $L^2(\R^n)$ with the function $G$ since $G$ is bounded and continuous.
In view of the Lemma below, one has
$ B(A - A')(Z , - \overline Z) = 0$, for all $Z\in\C^n$. One finishes the proof similarly to the one of Theorem  \ref{t1}. \fpr

  \begin{lemm} Let  $G$ be a bounded continuous function on $\R^{2n}$ and set
 $A' =  Op^{AW} (G) $ (bounded operator in $L^2(\R^n)$).
Then:
%----
\be\label{BSAW} B(A') (Z , - \overline Z) = (S_{1}G) (Z). \ee
%-----

 \end{lemm}

 {\it Proof of the Lemma.} From (\ref{egalite}) with $A$ replaced by $A'$,
one has, for any  $X\in\C^n$:
 %----
 $$   < A' \Psi_X  , \Psi_{-X}  >  = (2\pi )^{-n}  \int _{\R^{2n}}  G(Z)
e^{ - \frac {1} {4} \big ( |X - Z|^2 + |X +  Z|^2 \big ) +  i
\sigma (X , Z)       } dZ. $$
%------
From  (\ref{bisymbole}), (\ref{PSEC}) and (\ref{S-lambda-1}) with $\lambda = 1$,
one then obtains (\ref{BSAW}). \fpr

 \section{Iterated commutators.}\label{s5}

{\it Proof of (\ref{segal-EC}).}  We first note that  $\Psi_0 (u) =  \pi ^{-n/4} e^{ - \frac {1} {2} |u|^2}$, 
and that, as $\Phi_S(X)$
is the operator defined in (\ref{PhiS}), we have:
%---
$$ e^{\Phi_S(X) }  \Psi_0 =  e^{ \frac{1} {2} |X|^2 }  \Psi_X $$
%----
and:
%---
$$  e^{\Phi_S(Z) }   e^{\Phi_S(X) }  = e^{ - \frac{i} {2} \sigma (Z, X) }   e^{\Phi_S(Z+ X ) }.   $$
%----
Thus, one gets (\ref{segal-EC}). \fpr

{\it Proof of Theorem \ref{transfo-beals}.}  From (\ref{segal-EC}):
%---
 $$ \Big < e^{- \Phi_S(Z) }  A e^{\Phi_S(Z) }   \Psi_{X},
 \Psi_{X}   \Big > =  e^{  |Z|^2 - i\sigma (Z, X)    }  <  A \Psi_{X+Z } ,\Psi_{X-Z } >. $$
%---
With equality (\ref{PSEC}):
%-----
$$  <   \Psi_{X+Z } ,\Psi_{X-Z } >  =  e^{ - |Z|^2 + i\sigma (Z, X)    }. $$
%-------
These two equalities prove Theorem \ref{transfo-beals}. \fpr

{\it Proof of Proposition \ref{beals}.} One has $\Psi_X = W_X \Psi_0$
with
%----
\be\label{WX}  (W_X f) (u) = f(u - x)  e^{  i u \cdot \xi -  \frac {i} {2} x \cdot \xi}.  \ee
%-----
Thus:
%---
\be\label{prod-scal}  <  A \Psi_{X+Z}  , \Psi_{X-Z}  > = < W_{Z-X} A  W_{Z+X} \Psi_0 , \Psi_0>.\ee
%----
Denote by $V_1, \dots, V_n$ the canonical basis of $\R^n$. Set $a(V_j) = u_j + \partial _{u_j}$
and $a^{\star} (V_j) = u_j - \partial _{u_j}$. One checks that:
%----
$$ [ a(V_j) , W_{x , \xi}  ]   = (x_j + i \xi_j) W_{x , \xi} .  $$
%------
For all $X = (x , \xi)$, set $\varphi_j (X) = x_j + i \xi_j$.
Let  $A$ be a bounded operator in  $L^2(\R^n)$ such that, the right hand side of (\ref{norme-beals}) is well defined.  
Since $a(V_j) \Psi_0 = 0$:
%----
 \begin{align*} <  W_{Z-X} A W_{Z+X} \Psi_0 , a^{\star} (V_j) \Psi_0 >
&= < [ a(V_j) , W_{Z - X}  ] A W_{Z+X} \Psi_0 , \Psi_0 >  \\
&\ + < W_{Z-X} [ a(V_j) , A]W_{Z+ X}  \Psi_0 , \Psi_0 > \\
&\ +<  W_{Z-X} A [ a(V_j) , W_{Z + X}  ] \Psi_0 , \Psi_0 > \\%-----
& = 2 \varphi _j (Z)  <  W_{Z-X} A W_{Z+X} \Psi_0 ,\Psi_0 > +
< W_{Z-X} [ a(V_j) , A]W_{Z+ X}  \Psi_0 , \Psi_0 > .
 \end{align*}
%-------
By iteration and using (\ref{prod-scal}),  one gets, for all multi-indices $\alpha$:
%----
$$ \varphi (Z) ^{\alpha } <  A \Psi_{X+Z}  , \Psi_{X-Z}  > = \sum _{\beta + \gamma = \alpha }
c _{\alpha, \beta, \gamma } <  W_{Z-X} ({\rm ad}\, a(V) ) ^{\beta} A W_{Z+X} \Psi_0 , (a^{\star} (V))^{\gamma}  \Psi_0 >.  $$
%---
One  set above $ \varphi (Z) ^{\alpha } = \prod \varphi _j (Z) ^{\alpha _j} $
and similarly for $ (a^{\star} (V))^{\gamma}$. Thus, there exists $C_n$ such that:
%---
$$  |<  A \Psi_{X+Z}  , \Psi_{X-Z}  >| \  \sum _{ |\alpha | \leq 2n+1 } |\varphi (Z) ^{\alpha } |
\leq C_n  \sum _{ |\beta | \leq 2n+1 }\Vert ({\rm ad} \ a(V)) ^{\beta} A \Vert.$$
%----
One then gets  (\ref{norme-beals}) since the inverse of the left hand side belongs to $L^1(\R^{2n})$.
\fpr

\end{document}